\documentclass{amsart}

\usepackage{enumerate}
\usepackage{amsmath}
\usepackage{amssymb}
\usepackage{amsthm}

\usepackage{graphicx}
\usepackage{type1cm}
\usepackage{eso-pic}
\usepackage{color}
\makeatletter

\usepackage{math}
\usepackage{thm}
\usepackage{ref}

\newcommand{\R}{\ensuremath{\mathbb{R}}}

\newcommand{\D}{\ensuremath{\mathcal{D}}}
\newcommand{\N}{\ensuremath{\mathbb{N}}}

\newcommand{\Div}{\operatorname{div}}
\newcommand{\DIV}{\operatorname{Div}}

\newcommand{\B}{\mathcal{B}}
\newcommand{\E}{\mathbb{E}}
\newcommand{\F}{\mathbb{F}}
\newcommand{\G}{\mathbb{G}}

\newcommand{\pa}{\partial}
\newcommand{\grad}{\operatorname{grad}}

\newcommand{\rg}{\operatorname{rg}}
\newcommand{\id}{\operatorname{id}}

\renewcommand{\span}{\operatorname{span}}
\renewcommand{\phi}{\varphi}

\title
	[Asymptotic stability of local Helfrich minimizers]
	{Asymptotic stability of local Helfrich minimizers}

\author
	[Daniel Lengeler]
	{Daniel Lengeler}

\address
	{Department of Mathematics \newline\indent
	 Universit{\"a}t Regensburg --
	 Universit{\"a}tsstr.~31, 93053 Regensburg, Germany}
	
\email
	{daniel.lengeler@mathematik.uni-regensburg.de}

\keywords
	{Willmore energy, Canham-Helfrich energy, gradient flow, geometric flow, Willmore flow, Helfrich flow, Helfrich equation, Stokes system, linear elliptic system, fluid dynamics, biological membrane, lipid bilayer, well-posedness, stability}
	 
\thanks{I gratefully acknowledge support by DFG SPP 1506 "Transport Processes at Fluidic Interfaces".}		 

\subjclass
	[2010]
	{Primary: 35Q92; Secondary: 35B35, 35Q74, 35K25, 76D27}

\date
	{\today}

\begin{document}
\begin{abstract}
We show that local minimizers of the Canham-Helfrich energy are asymptotically stable with respect to a model for relaxational fluid vesicle dynamics that we already studied in previous papers (\cite{lengeler, lengelerwellp}). The proof is based on a \L ojasiewicz-Simon inequality.
\end{abstract}
\maketitle

\section*{Introduction}
In \cite{lengeler, lengelerwellp} we started the analysis of a basic model of fluid vesicle dynamics. In this model we describe the evolution of biological vesicles by considering a homogeneous, Newtonian surface fluid (\cite{scriven,aris}) subject to suitable elastic stresses that is immersed in a homogeneous, Newtonian bulk fluid. For a detailed introduction to the physics and mathematics of fluid vesicles we refer the reader to \cite{lengeler, lengelerwellp} and the references therein. There we showed that for most applications one can safely neglect inertial forces and, hence, restrict the model to purely relaxational dynamics. In this case the model takes the form
\begin{equation}\label{eqn:final}
 \begin{aligned}
 \Div S&=0&&\mbox{ in }\Omega\setminus\Gamma_t,\\
 \Div u&=0&&\mbox{ in }\Omega\setminus\Gamma_t,\\
 \DIV {^fT} + [\![S]\!]\nu_t&=-\DIV {^eT} &&\mbox{ on }\Gamma_t,\\
 \DIV u&=0&&\mbox{ on }\Gamma_t,\\
 u&=0&&\mbox{ on }\pa\Omega.
 \end{aligned}
 \end{equation}
Here, $\Omega$ is a smooth bounded domain in $\R^3$ containing a closed moving vesicle $\Gamma_t$, $\nu_t$ is the outer unit normal on $\Gamma_t$, and $u$ is the velocity of the bulk fluid in $\Omega\setminus\Gamma_t$ and the velocity of the vesicle on $\Gamma_t$; these velocities are assumed to coincide on $\Gamma_t$, that is, $u$ is a continuous function. Furthermore, $S=2\mu_b Du - \pi I$ is the Newtonian bulk stress tensor with the constant dynamic viscosity $\mu_b$ of the bulk fluid, the symmetric part $Du$ of the gradient of $u$, and the bulk pressure $\pi$, 
$[\![S]\!]$ is the jump of the bulk stress tensor across the membrane (subtracting the outer limit from the inner limit), $\DIV$ is the surface divergence (see below), and $T={^fT}+{^eT}$ is the surface stress tensor which is composed of a fluid part ${^fT}$ and an elastic part ${^eT}$. In coordinates we have $^fT^i_\alpha={^f\tilde T_\alpha^\beta}\pa_\beta^i$ with (see \cite{scriven,aris,lengeler,lengelerwellp})
\[^f\tilde T_\alpha^\beta=-q\,\delta_\alpha^\beta+2\mu\,(\D u)_\alpha^\beta=-q\,\delta_\alpha^\beta+\mu\,g^{\beta\gamma}(v_{\alpha;\gamma}+v_{\gamma;\alpha}-2w\,k_{\alpha\gamma})\]
and 
\[^{h}T^i_\alpha=\kappa\,\big((H-C_0)^2/2\,\pa_\alpha^i-(H-C_0)\,k_\alpha^\beta\pa_\beta^i-(H-C_0)_{,\alpha}\nu^i_t\big).\]
Here, $q$ is the surface pressure acting as a Lagrange multiplier with respect to the constraint $\DIV u=0$, $\mu$ is the constant dynamic viscosity of the surface fluid, $\D u$ is the \emph{surface rate-of-strain tensor}, $k$, $H$, and $K$ are the second fundamental form, twice the mean curvature, and the Gauss curvature of $\Gamma_t$, respectively, $\kappa$ is the bending rigidity, $C_0$ is the spontaneous curvature, $\partial_\alpha$ denotes the $\alpha$-th coordinate vector field, and the semicolon denotes covariant differentiation while the comma indicates usual partial differentiation. Furthermore, on $\Gamma_t$ we decomposed the function $u=v+w\,\nu_t$ into its tangential and its normal part. Throughout the paper, latin indices refer to Cartesian coordinates in $\R^3$ while greek indices refer to arbitrary coordinates on $\Gamma_t$. In particular, we note that the surface stress tensors are instances of \emph{hybrid tensor fields} (\cite{scriven,aris}) taking a tangential direction and returning a 
force density that is, in general, not tangential. The surface divergences for the non-tangential vector field $u$ and the hybrid tensor field $T$ 
can 
be written as
\begin{equation*}
\begin{aligned}
\DIV u &=g^{\alpha\beta}\langle\pa_\alpha u,\pa_\beta\rangle_{e},\\
(\DIV T)^i&=g^{\alpha\beta} T^i_{\alpha;\beta},
\end{aligned}
\end{equation*}
where $g$ denotes the Riemannian metric on $\Gamma_t$ induced by the Euclidean metric $e$ in $\R^3$, and the semicolon denotes the corresponding covariant differentiation of the covectors $(T^i_{\alpha})_{\alpha=1,2}$ (for fixed $i$). We showed in \cite{lengeler} that
\begin{equation}\label{fullt}
\begin{aligned}
\DIV u=&\Div_g v - w\, H,\\
 \DIV {^fT} = &-\grad_g q -q\,H\nu_t + \mu\,\big(\Delta_g v + \grad_g(w\,H) + K v -2\Div_g(w\,k)\big)\\ &  +2\mu\big(\langle\nabla^g v,k\rangle_g-w\,(H^2-2K)\big)\nu_t\\
\DIV {^eT}=& -\kappa\big(\Delta_g H + H(H^2/2-2K)+C_0(2K-HC_0/2)\big)\nu_t.
\end{aligned}
\end{equation}
Here, $\grad_g$, $\Div_g$, $\nabla^g$, $\Delta_g$ denote the differential operators (acting on tangential tensor fields) corresponding to the metric $g$, and, with a slight abuse of notation, we write $\langle\nabla^g v,k\rangle_g$ for the contraction of the tensor fields $\nabla^g v$ and $k$ using  $g$.

At first sight, the basic structure of our system \eqref{eqn:final} might not be so clear. Note that $\DIV {^eT}$ can be computed from $\Gamma_t$ alone. Hence, we have to solve the Stokes-type system defined by the left hand side of \eqref{eqn:final} with $-\DIV {^eT}$ as a right hand side for the fluid velocity $u$. Then, the normal part $w$ of $u$ on $\Gamma_t$ tells us how the vesicle will move in the next instant. Furthermore, it is easy to conclude from \eqref{eqn:final}$_{2,4}$ that area and the enclosed volume of each connected component $\Gamma^i_t$ of $\Gamma_t$ are preserved under this evolution; see \cite{lengeler}. Hence, the phase space $N$ of our system consists of embedded surfaces $\Gamma\subset\Omega$ with a fixed number $k$ of connected components $\Gamma^i$ and with fixed area and enclosed volume of each $\Gamma^i$. Moreover, we showed in \cite{lengeler} that \eqref{eqn:final} can be formulated as a gradient flow of the \emph{Canham-Helfrich energy} (see \cite{canham,helfrich,evans})
\begin{equation*}
 F(\Gamma)=\int_\Gamma \frac\kappa 2(H-C_0)^2\ dA
\end{equation*}
with respect to a suitably defined Riemannian metric on $N$. This observation will be particularly important for the present paper. For this reason we will repeat the details. For $\Gamma\in N$, the tangent space $T_\Gamma N$ can be identified with the space of scalar fields $w$ on $\Gamma$ such that the linearized constraints
\begin{equation}\label{eqn:tmp}
\int_{\Gamma^i} w\, dA=0\quad\text{ and }\quad\int_{\Gamma^i} w\,H\, dA=0
\end{equation}
hold for all $i=1,\ldots,k$. For $w\in T_\Gamma N$, consider the system
\begin{equation}\label{eqn:esystem}
\begin{aligned}
\Div S&=0 &&\mbox{ in }\Omega\setminus \Gamma,\\
\Div u&=0&&\mbox{ in }\Omega\setminus \Gamma,\\
P_\Gamma(\DIV {^fT} + [\![S]\!]\nu) &= 0&&\mbox{ on }\Gamma,\\
\DIV u&=0&&\mbox{ on }\Gamma,\\
u\cdot\nu&=w&&\mbox{ on }\Gamma,\\
u&=0&&\mbox{ on }\pa\Omega.
\end{aligned}
\end{equation}
Here, $P_\Gamma$ denotes the field of orthogonal projections onto the tangent spaces of $\Gamma$. Note that the conditions \eqref{eqn:tmp} are necessary for the solvability of these equations, due to the incompressibility constraints. For $w_1,w_2\in T_\Gamma N$, define the Riemannian metric on $N$ associated with fluid vesicle dynamics by
\begin{equation}\label{eqn:rmetric}
 \langle w_1,w_2\rangle_V:=2\mu_b\int_{\Omega\setminus\Gamma}\langle Du_1,Du_2\rangle_e\,dx+2\mu\int_{\Gamma}\langle \D u_1,\D u_2\rangle_g\,dA,
\end{equation}
where $u_1,u_2$ solve the system \eqref{eqn:esystem} with data $w_1,w_2$. Note that the length of a curve in $N$ endowed with this metric is given by the energy dissipated during the corresponding forced deformation of the membrane. The representation of $-dF_\Gamma$ with respect to the metric \eqref{eqn:rmetric} is given by $[u]_\Gamma\cdot\nu$, where $u$ solves \eqref{eqn:final} and where we use the notation $[u]_\Gamma$ to emphasize that we are taking the trace of $u$ on $\Gamma$.
Indeed, for all $w\in T_\Gamma N$ and corresponding solutions $\tilde u$ of \eqref{eqn:esystem} we have
\begin{equation}\label{gradient}
\begin{aligned}
\langle u\cdot\nu,w\rangle_V&=2\mu_b\int_{\Omega\setminus\Gamma}\langle Du,D\tilde u\rangle_e\,dx+2\mu\int_{\Gamma}\langle \D u,\D \tilde u\rangle_g\,dA\\
&=-\int_{\Gamma}\langle[\![S]\!]\nu,\tilde u\rangle_{e}\,dA-\int_{\Gamma}\langle\DIV {^fT},\tilde u\rangle_{e}\,dA=\int_{\Gamma}\langle\DIV {^hT},\tilde u\rangle_{e}\,dA\\
&=-\int_\Gamma \grad_{L_2}F_\Gamma\,w\,dA=-dF_\Gamma(w).
\end{aligned}
\end{equation}
Here, $S$ and $^fT$ denote the stress tensors with respect to $u$, and we used integration by parts for the second identity (see \cite{lengeler}), \eqref{eqn:final}$_3$ for the third identity, and \eqref{fullt}$_3$ for fourth identity.
We conclude that, indeed, \eqref{eqn:final} is the gradient flow of the Canham-Helfrich energy on $N$ endowed with the Riemannian metric \eqref{eqn:rmetric}. In particular, the energy $F$ is a strict Lyapunov functional, and, along the flow,
\begin{equation*}
\frac{d}{dt} F(\Gamma_t)=dF_{\Gamma_t}(u\cdot\nu)=-\langle u\cdot\nu,u\cdot\nu\rangle_V=-2\mu_b\int_{\Omega\setminus \Gamma_t}|Du|^2_{e}\,dx-2\mu\int_{\Gamma_t}|\D u|^2_g\,dA.
\end{equation*}
In \cite{lengeler} we concluded from this identity that the equilibria of \eqref{eqn:final} must satisfy the well-known Helfrich equation
\begin{equation}\label{helfrich}
 \grad_{L_2}F_\Gamma+[\![\pi]\!]+q\,H=0
\end{equation}
with the pressure jump and the surface pressure acting as Lagrange multipliers with respect to the volume and area constraints.

We refer the reader to the introduction in \cite{lengeler} for an overview of what is known concerning equilibrium configurations and dynamics of fluid vesicles from a physical point of view on the one hand and from the point of view of a rigorous mathematical analysis on the other hand. We only remark that in \cite{lengeler} we presented a thorough $L_2$-analysis of the Stokes-type systems defined by the left hand sides of \eqref{eqn:final} and \eqref{eqn:esystem} that takes into account geometric variations of the membrane, while in \cite{lengelerwellp} we proved a local well-posedness result for the dynamical system \eqref{eqn:final}; both results will turn out to be fundamental ingredients in the present article. Here, we will show that surfaces that minimize the Canham-Helfrich energy locally are asymptotically stable, that is, solutions starting near such a minimizer $\Gamma_0$ exist for all times,  remain nearby, and converge to a possibly different minimizer. The limit will in general be different from $\Gamma_0$ 
because the Canham-Helfrich energy admits continuous symmetries, and the equilibria therefore constitute a finite-dimensional manifold in phase space. In fact, these symmetries do not only include the rather trivial Euclidean group but, in the case of vanishing spontaneous curvature and higher genus, also \emph{special conformal transformations}; see \cite{seifert97} and the references therein.\footnote{These symmetries can in fact be observed under the microscope where thermal excitation continuously drives the vesicle along the manifold of equilibria; this phenomenon is called conformal diffusion.} For parabolic (that is, purely relaxational) systems, stability usually follows in a more or less straightforward manner from the well-known principle of linearized stability provided that the spectrum of the linearization is strictly negative; see for instance \cite{amann90,lunardi95}. However, due to the symmetries mentioned above, in our case the kernel of the linearization will be non-trivial; this is a 
typical situation in geometric problems. One can deal with this difficulty by center manifold theory (which is technical; see for instance \cite{henry81,lunardi95}) or by the \emph{generalized principle of linearized stability} (see \cite{pruess09}). For the latter, however, one needs to know quite a lot about the equilibria and the symmetry group; more precisely, one needs to assure that the linearization's kernel not only contains the tangent space to the manifold of equilibria (which is always true) but actually coincides with it. Usually, this is shown by direct computations. Proving this in our case is difficult for two reasons. Firstly, as mentioned above, the symmetry group can be rather complicated, and secondly, almost all local minimizers have a highly non-trivial configuration; in fact, while the round sphere is a solution of \eqref{helfrich} for any choice of $C_0$, it is the only known solution of spherical topology for $C_0=0$ which admits an analytical expression. Fortunately, for gradient-type 
systems there exists a third method for proving stability, the \L ojasiewicz inequality (\cite{loja62,loja84}) and its infinite-dimensional analogue, the \L ojasiewicz-Simon inequality (\cite{simon83,jendoubi98}). The coincidence of the tangent space to the manifold of equilibria and the linearization's kernel essentially says that transversely to the manifold of equilibria the energy grows quadratically; this leads, in fact, to exponential convergence. However, as long as the energy grows to some (uniformly bounded) power in transverse directions, one at least has algebraic convergence; this essentially is the content of the \L ojasiewicz inequality and its application to stability. In finite dimensions an analytic energy always grows to some power in transverse directions since otherwise it would be constant in this direction. Now, the essential step towards an application in infinite dimensions is to note that all the critical directions are contained in the linearization's kernel which in applications 
usually is finite-dimensional. An abstract result following such a Lyapunov-Schmidt type reduction which is sufficiently general for our purpose can be found in \cite{chill03}.

The present article is organized as follows. In section \ref{loja} we prove a suitable \L ojasiewicz-Simon inequality for our system by applying results from \cite{chill03}. In section \ref{asym} we will combine this inequality with the local well-posedness result from \cite{lengelerwellp} and a quantitative form of parabolic regularization to prove asymptotic stability. In the appendix we derive the second variation of the Canham-Helfrich energy and the linearization of its $L_2$-gradient.

Before we proceed, let us fix some notation. Throughout the article, let $\Omega\subset\R^3$  be a smooth bounded domain. For a closed surface $\Gamma\subset\Omega$ we write $\nu$ for outer unit normal, $\Gamma^i$, $i=1,\ldots,k$, for the connected components of $\Gamma$, $\Omega^i$ for the open set enclosed by $\Gamma^i$, and we define \[\Omega^0:=\Omega\setminus\big(\bigcup_{i=1}^k \Gamma^i\cup\Omega^i\big).\]
We denote by $[u]_\Gamma$ the trace of the bulk field $u$ on $\Gamma$; however, when there is no danger of confusion we will sometimes omit the brackets. Sometimes we use the notation $u\cdot v$ instead of $\langle u,v\rangle_e$ for $u,v\in\R^3$. We denote by $r(a)$ generic tensor fields that are polynomial or analytic functions of their argument $a$. For tensor fields $r_1$ and $r_2$ we write $r_1*r_2$ for any tensor field that depends in a bilinear way on $r_1$ and $r_2$, and we use the abbreviations $r*(r_1,\ldots,r_k)=r*r_1 + \ldots + r*r_k$ and $r^k=r*\ldots*r$ (with $k$ factors on the right hand side).
For $p\in (1,\infty)$, $k\in\N$, and $s\in\R_+\setminus\N$ we denote by $W^k_p$ the usual Sobolev spaces and by $W^s_p$ the Sobolev-Slobodetskij spaces. For an arbitrary smooth, $d$-dimensional Riemannian manifold $(M,\tilde e)$ the norm of the latter spaces is given by
\[\|T\|_{W^s_p(M)}=\|T\|_{W^k_p(M)}+|(\nabla^{\tilde e})^k T|_{W^{s-k}_p(M)},\]
where $k$ is the largest integer smaller than $s$ and 
\[|(\nabla^{\tilde e})^k T|_{W^{s-k}_p(M)}^p:=\int_M\int_M \frac{|(\nabla^{\tilde e})^k T(x)-(\nabla^{\tilde e})^k T(y)|^p_e}{d_{\tilde e}(x,y)^{d+(s-k)p}}\,dV_{\tilde e}(x)\,dV_{\tilde e}(y).\]
In this formula $d_{\tilde e}$ is the Riemannian distance function while $dV_{\tilde e}$ is the volume element corresponding to $\tilde e$. Finally, we define $H^s:=W^s_2$.

\section{\L ojasiewicz-Simon Inequality}\label{loja}
Let $M$ be the (formal) manifold of closed, embedded surfaces of class $H^{7/2}$ contained in $\Omega$. We introduce local coordinates near an arbitrary element of $M$ by approximating it by a smooth $\Gamma\in M$ and writing nearby surfaces as graphs over $\Gamma$. Let us make this more precise. We denote by $S_{\alpha}$, $\alpha>0$, the open set of points in $\Omega$ whose distance from $\Gamma$ is less than $\alpha$. It's a well-known fact from elementary differential geometry that there exists a maximal $\kappa_\Gamma>0$ such that the mapping
\begin{equation*}
 \begin{aligned}
\Lambda: \Gamma\times (-\kappa_\Gamma,\kappa_\Gamma)\rightarrow
S_{\kappa_\Gamma},\
(x,d)\mapsto x + d\,\nu(x)  
 \end{aligned}
\end{equation*}
is a diffeomorphism. 
For $h\in U=U(\Gamma):=\{h\in H^{7/2}(\Gamma)\,|\,|h|<\kappa_\Gamma/2\text{ on }\Gamma\}$, let 
\[\phi_h:\Gamma\rightarrow\R^3,\ x\mapsto x+h(x)\nu(x)\]
and $\Gamma_h:=\phi_h(\Gamma)\subset\R^3$. Then, the bijection
\[\Psi: U\rightarrow \Psi(V)\subset M,\quad h\mapsto \Gamma_h\]
defines a local chart. We consider the Canham-Helfrich energy $F$ as a function on $U$, that is, for $h\in U$, we let 
\[\tilde F(h):=(F\circ\Psi)(h)=\frac\kappa 2\int_{\Gamma}(H(h)-C_0)^2\,\sqrt{g(h)}\,dA.\]
Here and in the following, we denote by $H(h)$, $g_{\alpha\beta}(h)$, etc. the geometric quantities on $\Gamma_h$ pulled back by $\phi_h$ and $\sqrt{g(h)}:=\sqrt{\det (g_{\alpha\beta}(h))}$. Furthermore, for $i=1,\ldots,k$ we define $\Gamma_h^i:=\phi_h(\Gamma^i)$,
\[A_i(h):=\int_{\Gamma_h^i}\,dA_h,\quad\text{and}\quad V_i(h):=\frac13\int_{\Gamma_h^i}x\cdot\nu_h\,dA_h(x),\]
where $\nu_h$ is the outer unit normal and $dA_h$ is the area element on $\Gamma_h$. Note that, by the divergence theorem, $V_i(h)$ is nothing but the volume of the set enclosed by $\Gamma_h^i$.

\begin{lemma}\label{lemma:ana}
Let $X,Y$ be smooth tangent vector fields on $\Gamma$. Then, the following functions are well-defined and analytic:
\begin{itemize}
\item[(i)] $U\rightarrow H^{5/2}(\Gamma):$\quad $h\mapsto g(h)(X,Y)$ and $h\mapsto \sqrt{g(h)}$,
\item[(ii)] $U\rightarrow H^{3/2}(\Gamma):$\quad $h\mapsto k(h)(X,Y)$, $h\mapsto H(h)$, and $h\mapsto K(h)$,
\item[(iii)] $U\rightarrow \R:$\quad $h\mapsto A_i(h)$, $h\mapsto V_i(h)$ (for $i=1,\ldots,k$), and $h\mapsto \tilde F(h)$,
\item[(iv)] $U\rightarrow (H^{1/2}(\Gamma))':$\quad $h\mapsto d\tilde F_h$.
\end{itemize}
\end{lemma}
\proof We proved in appendix B in \cite{lengeler} that 
\begin{equation*}
\begin{aligned}
 g_{\alpha\beta}(h)&=r_{\alpha\beta}(h/\kappa_\Gamma,hk,\nabla h),\\
 g^{\alpha\beta}(h)&=r^{\alpha\beta}(h/\kappa_\Gamma,hk,\nabla h),\\
 k_{\alpha\beta}(h)&=r_{\alpha\beta}(h/\kappa_\Gamma,hk,\nabla h)*\big(k,\nabla h/\kappa_\Gamma,h\nabla^g k,(\nabla^g)^2 h\big).
\end{aligned}
\end{equation*}
Since $\det(g_{\alpha\beta}(h))$ is uniformly positive and $H^s(\Gamma)$ is an algebra for $s>1$, claim (i) follows. From these considerations and $H(h)=g^{\alpha\beta}(h)\,k_{\alpha\beta}(h)$, $K(h)=\det(k_{\alpha\beta}(h))$ we infer that claim (ii) is true. Since the map
\[L^1(\Gamma) \rightarrow \R,\quad f\mapsto \int_{\Gamma}f\,dA\]
is linear and bounded, from (i) and (ii) we can deduce the first and the third claim in (iii). Furthermore, by a straightforward computation one can check that the normal $\nu_h$ to $\Gamma_h$ is given by
\begin{equation}\label{nuh}
\nu_h\circ\varphi_h=\frac{\nu-\grad_{\bar g}h}{(1+|\grad_{\bar g}h|_g^2)^{1/2}}=r(hk,\nabla h) 
\end{equation}
with the metric $\bar g_{\alpha\beta}=g_{\alpha\beta}-2hk_{\alpha\beta}+h^2k_{\alpha\gamma}k^\gamma_{\beta}$. In view of
\[V_i(h)=\frac13\int_{\Gamma_h^i}x\cdot\nu_h\,dA_h=\frac13\int_{\Gamma^i}(y+h\,\nu)\cdot(\nu_h\circ\phi_h)\,\sqrt{g(h)}\,dA(y),\]
this proves the second claim in (iii). Finally, note that for $w\in H^{7/2}(\Gamma)$ we have
\[d\Psi_h(w)=\nu_h\cdot\big((w\nu)\circ\phi^{-1}_h\big)\]
and thus for sufficiently smooth $h\in U$
\begin{equation}\label{Fgrad}
\begin{aligned}
d\tilde F_h(w)&=dF_{\Gamma_h}d\Psi_h(w)\\
&=\kappa\int_{\Gamma}\big(\Delta_{g(h)} H(h) + H(h)((H(h)^2/2-2K(h))\\
&\quad\quad\quad+ C_0(2K(h)-H(h)\,C_0/2))\big)(\nu_h\circ\phi_h)\cdot\nu\,w\,\sqrt{g(h)}\,dA.
\end{aligned}
\end{equation}
In view of (i), (ii), \eqref{nuh}, and the fact that $H^{3/2}(\Gamma)$ is an algebra, we see that the summands in \eqref{Fgrad} not involving the Laplacian give rise to an analytic map from $U$ to $H^{3/2}(\Gamma)$. We proved in appendix B in \cite{lengeler} that
\begin{equation}\label{duri}
\begin{aligned}
\Delta_{g(h)}H(h)&=r(h/\kappa_\Gamma,hk,\nabla h)*(\nabla^g)^2 H(h) \\
&\quad+ r(h/\kappa_\Gamma,hk,\nabla h)*\big(k,\nabla h/\kappa_\Gamma,h\nabla^g k,(\nabla^g)^2 h\big)*\nabla^gH(h).
\end{aligned}
\end{equation}
Since
\[H^{3/2}(\Gamma)\cdot H^{1/2}(\Gamma)\hookrightarrow H^{1/2}(\Gamma),\]
that is, pointwise multiplication is continuous in the indicated function spaces,
the expressions in \eqref{Fgrad} involving the second summand in \eqref{duri} give rise to an analytic map from $U$ to $H^{1/2}(\Gamma)$. Finally, concerning the expressions involving the first summand in \eqref{duri}, we note that, on the one hand, since
\[H^{5/2}(\Gamma)\cdot L^{2}(\Gamma)\hookrightarrow L^{2}(\Gamma),\]
the map
\[H^{7/2}(\Gamma)\rightarrow \mathcal{B}(H^{2}(\Gamma),L^{2}(\Gamma)), \quad h\mapsto r(h/\kappa_\Gamma,hk,\nabla h)*(\nabla^g)^2\]
is analytic. On the other hand, for scalar fields $f$ integration by parts gives
\begin{equation*}
\begin{aligned}
\int_{\Gamma}r(h/\kappa_\Gamma,hk,&\nabla h)*(\nabla^g)^2f\,w\,dA\\
&=\int_{\Gamma}r(h/\kappa_\Gamma,hk,\nabla h)*\nabla f*\nabla w\,dA\\
&\quad+ \int_\Gamma r(h/\kappa_\Gamma,hk,\nabla h)*\big(k,\nabla h/\kappa_\Gamma,h\nabla^g k,(\nabla^g)^2 h\big)*\nabla f\,w\,dA,
\end{aligned}
\end{equation*}
proving that
\[H^{7/2}(\Gamma)\rightarrow \mathcal{B}(H^{1}(\Gamma),(H^{1}(\Gamma))'), \quad h\mapsto r(h/\kappa_\Gamma,hk,\nabla h)*(\nabla^g)^2\]
is also analytic. Thus, by interpolation, this map is analytic from $H^{7/2}(\Gamma)$ to $\mathcal{B}(H^{3/2}(\Gamma),(H^{1/2}(\Gamma))')$. This together with the second assertion in (ii) proves (iv).\qed

\begin{remark}\label{psi}
 The preceding proof shows that for $h\in U$ the (formal) derivative of $d\Psi_h$ defines linear isomorphisms from $H^s(\Gamma)$ to $H^s(\Gamma_h)$ for all $s\in [0,5/2]$ which are uniformly bounded in both directions for $h\in U$ being uniformly bounded in $H^{7/2}(\Gamma)$.
\end{remark}

For $A=(A_1,\ldots,A_k)$, and $V=(V_1,\ldots,V_k)$ with $A_i,V_i>0$, let $N=N_{A,V}\subset M$ denote the (formal) submanifold of closed, embedded surfaces of class $H^{7/2}$ contained in $\Omega$ which consist of $k$ connected components $\Gamma^i$ of fixed area $A_i$ and fixed enclosed volume $V_i$. For fixed smooth $\Gamma\in N$, in a neighbourhood of $\Gamma$ we can consider $N$ as an analytic submanifold of $U(\Gamma)$. By Lemma \ref{lemma:ana} (iii), the map \[f:U\rightarrow \mathcal{R}_1\times\ldots\times\mathcal{R}_k,\quad h\mapsto (f_1(h),\ldots,f_k(h)),\]
where $f_i(h)=(A_i(h),V_i(h))$, $\mathcal{R}_i=\R^2$ if $\Gamma^i$ is not a round sphere and $f_i(h)=A_i(h)$, $\mathcal{R}_i=\R$ else, is analytic, and its differential at $h=0$ is surjective. Indeed, for $w\in H^{7/2}(\Gamma)$ we have 
 \[d(f_i)_{h=0}(w)=\Big(\int_{\Gamma^i}w\,H\, dA,\int_{\Gamma^i} w\,dA\Big)\]
if $\Gamma_i$ is not a round sphere, and else
\[d(f_i)_{h=0}(w)=\int_{\Gamma_i} w\,dA;\]
 note that in the first case the functions $1$ and $H$ are linearly independent. Thus, by the implicit function theorem (see \cite{deimling}, for instance), there exist a closed complement $B$ of $T_\Gamma N:=\ker df_{h=0}$ in $H^{7/2}(\Gamma)$, bounded open neighbourhoods $\tilde U$ and $\hat U$ of the origin in $T_\Gamma N$ and $B$, respectively, and an analytic function $\gamma:\tilde U\rightarrow \hat U$ such that $\gamma(0)=0$, $d\gamma(0)=0$, and \[\{h\in \tilde U\times \hat U\,|\,f(h)=f(0)\}=\{w+\gamma(w)\,|\,w\in \tilde U\}=:G.\]
 Concatenating the map $1+\gamma:\tilde U\rightarrow G$, which is bianalytic, with the map $\Psi:U\rightarrow M$ yields local coordinates for $N$.
 
  For $\Gamma\in N$ and $s\ge 0$ let 
\[H_{n}^s(\Gamma):=\{w\in H^{s}(\Gamma)\ \text{such that \eqref{eqn:tmp} holds}\}.\]
Furthermore, let $L^2_n(\Gamma):=H_n^0(\Gamma)$, $Y_\Gamma:=H^{1/2}_n(\Gamma)$, and note that $T_\Gamma N=H^{7/2}_n(\Gamma)$. 

 \begin{remark}\label{gamma}
 For $w_0\in\tilde U$, let us have a closer look at the map
 \[d(1+\gamma)_{w_0}:T_\Gamma N\rightarrow T_{(1+\gamma)(w_0)}G,\]
 which is a linear isomorphism. Let $b_1,\ldots,b_n$ be a basis for $B$. From the definition of the functional $f$ it is not hard to see that there exist analytic maps $g_i:U\rightarrow H^{3/2}$, $i=1,\ldots,n$, such that
 \[T_{(1+\gamma)(w_0)}G=\ker df_{(1+\gamma)(w_0)}=Z_{w_0}\cap H^{7/2}(\Gamma)\]
 for the $L^2$-orthogonal complement $Z_{w_0}=(\span\{g_1(w_0),\ldots,g_n(w_0)\})^{\perp}\subset L^2(\Gamma)$. Thus, we have 
 \[d(1+\gamma)_{w_0}:w\mapsto w+a^i(w)\,b_i,\]
 where the real coefficients $a^i$ are determined by the linear system \[a^i(w)\,(b_i,g_j(w_0))_{L^2(\Gamma)}=-(w,g_j(w_0))_{L^2(\Gamma)}\] for $j=1,\ldots,n$. The matrix on the left hand side is invertible, all coefficients depend analytically on $w_0$, and the right hand side as a function in $w$ is a continuous functional on $L^1(\Gamma)$. In particular, for all $s\ge 0$, $d(1+\gamma)$ can be extended to an analytic map
 \[d(1+\gamma):\tilde U\rightarrow \B(H^{s}(\Gamma),H^{s}(\Gamma)),\]
 and for all $w_0\in \tilde U$
 \[d(1+\gamma)_{w_0}:H^s_{n}(\Gamma)\rightarrow Z_{w_0}\cap H^{s}(\Gamma)\]
 is an isomorphism with uniformly bounded continuity constants.
Furthermore, from Remark \ref{psi} we infer that $d\Psi_{(1+\gamma)(w_0)}: Z_{w_0}\cap H^s(\Gamma)\rightarrow H^s_n(\Gamma_{(1+\gamma)(w_0)})$, $s\in [0,5/2]$, defines isomorphisms which are uniformly bounded in both directions for $w_0\in \tilde U$.
 \end{remark}
 
Now we are ready to state one of the main results of this section.

\begin{theorem}[\L ojasiewicz-Simon inequality]\label{thm:ls}
 Let $\Gamma_0\in N$ be a smooth stationary point for $F$ in $N$, that is, $dF_{\Gamma_0}=0$. Then there exists a neighbourhood $\tilde U$ of $\Gamma_0$ in $N$ (in the above sense), a constant $c>0$, and a number $\theta\in (0,1/2]$ such that
 \begin{equation*}
  |F(\Gamma)-F(\Gamma_0)|^{1-\theta}\le c\|dF_\Gamma\|_{Y'_\Gamma}
 \end{equation*}
for all $\Gamma\in \tilde U$.
\end{theorem}

In order to prove the theorem, we have to analyze the second variation of $F$ in $N$. 
\begin{proposition}\label{fred}
For smooth $\Gamma\in N$ the second variation $d^2F_\Gamma:T_\Gamma N\rightarrow Y_{\Gamma}'$ is a Fredholm operator of index $0$.\footnote{The condition of smoothness can be weaked of course. It is assumed here only for simplicity.} In particular, we have 
 \begin{equation}\label{dim}
  \rg d^2F_\Gamma= (\ker d^2F_\Gamma)^\perp \cap Y_{\Gamma}'.
 \end{equation}
\end{proposition}
\proof In the appendix we show that the second variation has the form
\begin{equation*}
d^2F_\Gamma(w,\tilde w)=\kappa\int_{\Gamma}\big(\Delta_g w\,\Delta_g\tilde w + Bw\, \tilde w\big)\,dA,
\end{equation*}
where $w,\tilde w\in T_\Gamma N$ and $Bw=(a^{\alpha\beta} w_{,\beta})_{;\alpha} + b\,w$.
The operator $B$ obviously maps $T_\Gamma N$ compactly into $Y_\Gamma'$.
Hence, it suffices to consider the biharmonic operator alone. By standard arguments based on $L_2$-theory for the Laplacian on $\Gamma$ and Riesz' representation theorem, for some $\eta>0$ the operator $\Delta_g^2+\eta:H^2_n(\Gamma)\rightarrow (H^2_n(\Gamma))'$ is an isomorphism. From this and again $L_2$-theory for the Laplacian, we have that $\Delta_g^2+\eta:H^4_n(\Gamma)\rightarrow (L^2_n(\Gamma))'$ is an isomorphism, too. Indeed, each element of $(L^2_n(\Gamma))'$ is an element of $(H^2_n(\Gamma))'$ and hence the image of some $H^2_n(\Gamma)$ function. However, $(\Delta_g^2+\eta)
u\in (L^2_n(\Gamma))'$ means that $(\Delta_g^2+\eta)u\in L^2_n(\Gamma)$ modulo some linear combination of the functions $1$ and $H$ on each connected component of $\Gamma$, thus, $\Delta_g^2u\in L^2(\Gamma)$. Now, $L_2$-theory tells us that, in fact, $u\in H^4_n(\Gamma)$, and hence, $\Delta_g^2+\eta:H^4_n(\Gamma)\rightarrow (L^2_n(\Gamma))'$ is surjective. On the other hand, injectivity is obvious. Interpolating these results we obtain that $\Delta_g^2+\eta:T_\Gamma N\rightarrow Y_\Gamma'$ is an isomorphism. Thus, by Fredholm's alternative, for $K:=\eta(\Delta_g^2+\eta)^{-1}$ we have that $\id-K$ is a Fredholm operator of index $0$ on $T_\Gamma N$. However, it is easy to see that $\ker(\id-K)=\ker\Delta_g^2\subset T_\Gamma N$ and $(\Delta_g^2+\eta)(\rg (\id-K))=\rg \Delta_g^2\subset Y_\Gamma'$. Hence, $\Delta_g^2:TN_\Gamma\rightarrow Y_\Gamma'$ is Fredholm of index $0$, and, thus, the same is true for $d^2F_\Gamma$.

Finally, the symmetry of $d^2F_\Gamma$ implies that 
\begin{equation}\label{ker}
 \rg d^2F_\Gamma\subset (\ker d^2F_\Gamma)^\perp \cap Y_{\Gamma}'.
\end{equation}
But, by
\[\ker d^2F_\Gamma\cap (\ker d^2F_\Gamma)^\perp=\{0\}\quad\text{in}\ Y'_\Gamma\]
and the vanishing Fredholm index, the inclusion in \eqref{ker} must be an identity, that is,  \eqref{dim} holds.\qed

\proof (of Theorem \ref{thm:ls}) We consider $F$ as a function on $\tilde U$, that is, for $w\in \tilde U$, we let $\hat F(w):=\tilde F(w+\gamma(w))$. By Lemma \ref{lemma:ana} (iii), $\hat F$ is analytic. Furthermore, for $w_0\in\tilde U$ and $w_1\in T_\Gamma N$ we have
\begin{equation}\label{deriv}
\begin{aligned}
d\hat F_{w_0}(w_1)&=d\tilde F_{(1+\gamma)(w_0)}d(1+\gamma)_{w_0}(w_1)\\
&=dF_{\Gamma_{(1+\gamma)(w_0)}}d\Psi_{(1+\gamma)(w_0)}d(1+\gamma)_{w_0}(w_1). 
\end{aligned}
\end{equation}
Thus, by Remark \ref{gamma} it is sufficient to show that there exists a neighbourhood $\tilde U\subset T_{\Gamma_0} N$ of the origin and a number $\theta\in (0,1/2]$ such that 
\[|\hat F(w_0)-\hat F(0))|^{1-\theta}\le c\,\|d\hat F_{w_0}\|_{Y_\Gamma'}\]
for all $w_0\in \tilde U$. 

We want to apply Corollary 3.11 in \cite{chill03} with $X=V=T_\Gamma N$ and $Y=W=Y_\Gamma'$.
From the first equality in \eqref{deriv}, Remark \ref{gamma}, and Lemma \ref{lemma:ana} (iv), it follows that the map $\tilde U\rightarrow Y_\Gamma', w_0\mapsto d\tilde F_{w_0}$ is analytic. Moreover, since $\Gamma$ is a stationary point for $F$ in $N$ and $d\gamma(0)=0$, for $w_0,w_1\in T_\Gamma N$ we have
\[d^2\hat F_{w=0}(w_0,w_1)=d^2 F_\Gamma(w_0,w_1).\]
Thus, by Proposition \ref{fred}, $\ker d^2 \hat F_{w=0}$ is finite-dimensional and
\begin{equation}\label{sec}
\rg d^2 \hat F_{w=0}= (\ker d^2 \hat F_{w=0})^\perp \cap Y_\Gamma'.
\end{equation}
Let $\tilde P:Y_\Gamma \rightarrow \ker d^2 \hat F_{w=0}$ be a continuous projection and consider its restriction $P$ to $T_\Gamma N$, that is, the projection $P:T_\Gamma N\rightarrow \ker d^2 \hat F_{w=0}$ along $\ker P\cap T_\Gamma N$. For $w\in T_\Gamma N$ and $y'\in Y_\Gamma'$ we have
\[\langle P'y',w\rangle_{Y_\Gamma}=\langle y',\tilde Pw\rangle_{Y_\Gamma}\le\|y'\|_{Y_\Gamma'}\|\tilde Pw\|_{Y_\Gamma}\le c\,\|y'\|_{Y_\Gamma '}\|w\|_{Y_\Gamma},\]
which, by denseness of $T_\Gamma N$ in $Y_\Gamma$, proves that $P'$ leaves $Y_\Gamma '$ invariant. Since $\ker P'=(\rg P)^\perp$, equation \eqref{sec} takes the form
\begin{equation*}
\rg d^2 \hat F_{w=0}= \ker P' \cap Y_\Gamma'.
\end{equation*}
Hence, we checked all assumptions in Corollary 3.11 in \cite{chill03}. This finishes the proof.
\qed

\section{Asymptotic stability}\label{asym}
For $T>0$ we define $I:=(0,T)$. In \cite{lengelerwellp} we proved the following theorem with the exception of \eqref{sigmaest} which expresses parabolic regularization in a quantitative form.

\begin{theorem}\label{thm}
Let $\Gamma\subset\Omega$ be a smooth, closed surface that contains no round sphere. For all $p\in (3,\infty)\setminus\{4\}$ and sufficiently small $\epsilon>0$ there exists a time $T>0$ such that for all height functions $h_0\in \bar B_\epsilon(0)\subset W^{5-4/p}_p(\Gamma)$ there exists a unique
 \begin{equation*}
 h\in W^1_p(I,\,W^{2 - 1/p}_p(\Gamma)) \cap L_p(I,\,W^{5 - 1/p}_p(\Gamma))
\end{equation*}
with $\|h\|_{L_\infty(I\times\Gamma)}\le\kappa_\Gamma/2$ as well as suitable hydrodynamic fields $(u,\pi,q)$ such that $\Gamma_t=\Gamma_{h(t)}$ and $(u,\pi,q)$ solve \eqref{eqn:final} in the time interval $I$ with initial value $\Gamma_{h_0}$. The map
 \begin{equation*}
\bar B_\epsilon\subset W^{5-4/p}_p(\Gamma)\rightarrow W^1_p(J,\,W^{2 - 1/p}_p(\Gamma)) \cap L_p(J,\,W^{5 - 1/p}_p(\Gamma)),\quad h_0\mapsto h
\end{equation*}
is Lipschitz continuous. Furthermore, for sufficiently small $\delta>0$ and all $t'\in (0,T)$ there exists a constant $c>0$ such that
 \begin{equation}\label{sigmaest}
\|h\|_{C([t',T],W^{5-4/p+\delta}_p(\Gamma))}\le c.
\end{equation}
\end{theorem}
For details concerning the hydrodynamic fields see \cite{lengelerwellp}; alternatively, these can be constructed (in the $L_2$-scale) by applying Theorem 3.6 or Theorem 3.7 from \cite{lengeler}.

If $\Gamma$ consists only of round spheres, then  \eqref{eqn:final} is uniquely solved by the constant-in-time solution with $\Gamma_t=\Gamma$, $u=0$, and suitably chosen pressure functions. So far, we cannot prove the well-posedness of our system in the vicinity of a round sphere and in the case that $\Gamma$ contains both round spheres and non-spheres. The reason for this lies in the different degrees of gauge freedom of the pressure functions for round spheres on the one hand and non-spheres on the other hand; see \cite{lengelerwellp}.

\proof In order to show \eqref{sigmaest} we have to repeat parts the proof of Theorem 3.1 in \cite{lengelerwellp} in a time-weighted setting. We will only briefly sketch the procedure. Without further explanation we will use the notation from \cite{lengelerwellp}. Let $p\in [2,\infty)$ and $\mu\in (1/p,1)$ such that \[2+3\mu-4/p\in (2-1/p,5-4/p)\] is not a natural number. 

First, we prove that the linearization
\begin{equation}\label{eqn:tildesystemlinear}
\begin{aligned}
\mu_b\Delta u-\grad\pi&=f_1&&\mbox{ in }\Omega\setminus\Gamma,\\
\Div u&=f_2&&\mbox{ in }\Omega\setminus\Gamma,\\
\mu\big(\Delta_{g} v + \grad_{g}(w\,H) + K\, v -2\Div_{g}(w\,k)\big)&\\
-\grad_{g} q+ P_\Gamma [\![S]\!]\nu&=f_3^\top&&\mbox{ on } \Gamma,\\
2\mu\big(\langle\nabla^{g} v,k\rangle_{g}
-w\,(H^2-2K)\big)- q\,H+[\![S]\!]\nu\cdot\nu-Ah&=f_3^\perp&&\mbox{ on }\Gamma,\\
\Div_{g} v- w\, H&=f_4&&\mbox{ on }\Gamma,\\
u-v-w\,\nu&=f_5&&\mbox{ on }\Gamma,\\
\partial_t h-w&=f_6&&\mbox{ on }\Gamma
\end{aligned}
\end{equation}
admits a unique solution $(u,v,w,\pi,q,h)$ with 
\begin{equation*}
	\begin{array}{c}
		t^{1-\mu}u \in L_p(I,W^2_p(\Omega\setminus\Gamma,\,\bR^3)\cap W^1_p(\Omega,\R^3)), \quad t^{1-\mu}v\in L_p(I,W^2_p(\Gamma,T\Gamma)), \\
		t^{1-\mu}w\in L_p(I,W^{2-1/p}_p(\Gamma)),\quad t^{1-\mu}\pi \in L_p(I,W^1_p(\Omega\setminus\Gamma)),\quad  t^{1-\mu}q \in L_p(I,W^1_p(\Gamma)),\\
t^{1-\mu}h\in L_p(I,W^{5-1/p}_p(\Gamma)),\quad t^{1-\mu}\pa_t h\in L_p(I,W^{2-1/p}_p(\Gamma))
	\end{array}
\end{equation*}
and $\int_\Omega \pi\,dx=0$
provided that the data $(f_1,\ldots,f_6,h_0)$ satisfy
\begin{equation*}
	\begin{array}{c}
		t^{1-\mu}f_1\in L_p(I,L_p(\Omega\setminus\Gamma,\R^3)),\quad t^{1-\mu}f_2\in L_p(I,W^1_p(\Omega\setminus\Gamma)),\\
		t^{1-\mu}f_3^\top\in L_p(I,L_p(\Gamma,T\Gamma)), \quad t^{1-\mu}f_3^\perp\in L_p(I,W^{1-1/p}_p(\Gamma)),\\
t^{1-\mu}f_4\in L_p(I,W^{1}_p(\Gamma)),\quad t^{1-\mu}f_5,\,t^{1-\mu}f_6\in L_p(I,W^{2-1/p}_p(\Gamma,\R^3)),\\
h_0\in W^{2+3\mu-4/p}_p(\Gamma)
	\end{array}
\end{equation*}
with $\int_\Omega f_2\, dx=0$. To this end, we eliminate $h_0$ by choosing an extension $\tilde h$ such that $t^{1-\mu}\tilde h\in L_p(I,W^{5-1/p}_p(\Gamma))$ and $t^{1-\mu}\pa_t \tilde h\in L_p(I,W^{2-1/p}_p(\Gamma))$; see Proposition 3.1 in \cite{pruess04}. Next, we eliminate $(f_1,\ldots,f_5)$ by solving the stationary system
\begin{equation}\label{eqn:tildesystemlinearstat}
\begin{aligned}
\mu_b\Delta u-\grad\pi&=f_1&&\mbox{ in }\Omega\setminus\Gamma,\\
\Div u&=f_2&&\mbox{ in }\Omega\setminus\Gamma,\\
\mu\big(\Delta_{g} v + \grad_{g}(w\,H) + K\, v -2\Div_{g}(w\,k)\big)&\\
-\grad_{g} q+ P_\Gamma [\![S]\!]\nu&=f_3^\top&&\mbox{ on } \Gamma,\\
2\mu\big(\langle\nabla^{g} v,k\rangle_{g}
-w\,(H^2-2K)\big)- q\,H+[\![S]\!]\nu\cdot\nu&=f_3^\perp&&\mbox{ on }\Gamma,\\
\Div_{g} v- w\, H&=f_4&&\mbox{ on }\Gamma,\\
u-v-w\,\nu&=f_5&&\mbox{ on }\Gamma
\end{aligned}
\end{equation}
for almost all $t\in I$. In the proof of Theorem 3.2 in \cite{lengelerwellp} we showed unique $L_p$-solvability of the principal linearization of this system in the double half-space; combining this   result with standard localization and transformation techniques as well as the $L_2$-theory proved in \cite{lengeler}, we easily obtain unique $L_p$-solvability for \eqref{eqn:tildesystemlinearstat}. Thus, it is sufficient to solve \eqref{eqn:tildesystemlinear} for $h$ with all data vanishing except for $f_6$. We can write this system in the form 
\begin{equation}\label{ndequ}
\pa_t h +  L h=f_6,\quad h(0)=0.
\end{equation}
Here, $L: D(L)\rightarrow X$ with $X=W^{2-1/p}_p(\Gamma)$ and $D(L)=W^{5-1/p}_p(\Gamma)$ is the linear operator that maps $h\in D(L)$ to $w=[u]_\Gamma\cdot\nu\in X$, where $u$ solves \eqref{eqn:tildesystemlinearstat} with $f_3^\perp=Ah$ and all other data vanishing. This operator is closed as can be seen from the $L_p$-theory for \eqref{eqn:tildesystemlinearstat} and for the Laplacian on $\Gamma$. Furthermore, we proved in \cite{lengeler} that for $f_6\in L_p(I,X)$ equation \eqref{ndequ} admits a unique solution $h\in L_p(I,D(L))$ with $\pa_t h\in L_p(I,X)$. From this, however, using a summation argument one easily deduces that the same assertion holds with $I=(0,\infty)$ if we replace $L$ by $L+\lambda$ for a sufficiently large $\lambda>0$. Now, from Theorem 2.4 in \cite{pruess04} we finally obtain the existence of a unique solution $h$ of \eqref{ndequ} such that $t^{1-\mu}h \in L_p(I,D(L))$ and $t^{1-\mu}\pa_t h \in L_p(I,X)$ provided that $t^{1-\mu}f_6 \in L_p(I,X)$.

Next, we have to repeat the contraction mapping argument from section 4 in \cite{lengelerwellp} in a time-weighted setting. Essentially, throughout the proof we simply replace the spaces $\E_p(T)$, $\G_p(T)$, and $\F_p(T)$ by the corresponding time-weighted spaces, which we denote by $\E_{p,\mu}(T)$, etc., and correspondingly the time-trace space $W^{5-4/p}_p(\Gamma)$ by $W^{2+3\mu-4/p}_p(\Gamma)$. However, we have to prove that for $p>3$ and $\mu>3/p$ the statement analogous to Lemma 4.1 in \cite{lengelerwellp} holds. The analogue of assertion (i) in Lemma 4.1 follows from Proposition 3.1 in \cite{pruess04} and the existence of a bounded extension operator $\E_{p,\mu}(T)\rightarrow \E_{p,\mu}(\infty)$; the latter can be constructed quite simply via reflection. The analogue of assertion (ii) in Lemma 4.1 follows from Theorem 4.2 in \cite{meyries}. Furthermore, noting that for $\mu \ge 1/p+2/3$ we have $W^{-1+3\mu-4/p}_p(\Gamma)\hookrightarrow W^{1-1/p}_p(\Gamma)$, the assertion analogous to Lemma 4.2 in \cite{lengelerwellp} can be shown. Now, we can follow the proof of Theorem 1.1 in \cite{lengelerwellp} line by line to see that the first part of Theorem \ref{thm} holds with $W^{5-4/p}_p(\Gamma)$ replaced by $W^{2+3\mu-4/p}_p(\Gamma)$ and $W^1_p(I,\,W^{2 - 1/p}_p(\Gamma)) \cap L_p(I,\,W^{5 - 1/p}_p(\Gamma))$ replaced by the corresponding time-weighted space.

Finally, let us fix a $p\in (3,\infty)\setminus\{4\}$. We choose $\tilde p\in (p,\infty)$ and $\mu\in(1/\tilde p,1)$ such that the well-posedness result we just proved holds and such that the spatial regularities of the time-trace spaces of $\E_p(T)$ and $\E_{\tilde p,\mu}(T)$ coincide, that is, $ 2+3\mu-4/\tilde p=5-4/p$. Thus, we can apply the well-posedness theorem in the time-weighted setting to all $h_0\in\bar B_\epsilon\subset W^{5-4/p}_p(\Gamma)$ for sufficiently small $\epsilon$ and obtain a solution $h\in \E_p(T)$ such that
\[\|h\|_{\E_{\tilde p,\mu}(T)}\le c\]
for some constant $c>0$. However, for all $t'\in (0,T)$ the space $\E_{\tilde p,\mu}(T)$ obviously embeds into $C([t',T],W^{5-4/\tilde p}_{\tilde p}(\Gamma))$; this proves \eqref{sigmaest}.
\qed
\medskip
%

Now, we can prove the asymptotic stability of local Helfrich minimizers.  
\begin{theorem}
Let $\Gamma\in N$ be a smooth local minimizer for $F$ in $N$ that contains no round sphere, and let $\tilde U\subset T_\Gamma N$ be the set from Theorem \ref{thm:ls}.  For all $\sigma>0$ there exists an $\epsilon>0$ such that for all \[h_0\in\big(\bar B_\epsilon(0)\cap (1+\gamma)(\tilde U)\big)\subset W^{5-4/p}_p(\Gamma)\] the corresponding solution $h$ from Theorem \ref{thm} exists for all times and satisfies $\|h(t)\|_{W^{5-4/p}_p(\Gamma)}\le \sigma$ as well as 
\begin{equation}\label{konvergenz}
\|h(t)-h_\infty\|_{W^{5-4/p}_p(\Gamma)}\le c\,t^{-\beta}
\end{equation}
for all $t>0$, constants $c, \beta>0$, and some $h_\infty\in (1+\gamma)(\tilde U)$ with $F(\Gamma_{h_\infty})=F(\Gamma)$.
\end{theorem}
 
\proof The proof proceeds in two steps. First, we use Theorem \ref{thm:ls} to show that if our solution exists for all times then its energy will converge to the local minimum and it will satifisfy an arbitrarily small bound in a low norm. Then, combining these insights and parabolic regularization in a bootstrap argument we can prove that the solution will indeed exist for all times and converge in phase space.

By the Lipschitz continuity of the solution map from Theorem \ref{thm}, we can choose $\epsilon$ so small that for all $h_0$ as in the assertion of the present theorem the corresponding solution $h$ satisfies $\|h(t)\|_{W^{5-4/p}_p(\Gamma)}\le \sigma$ for all $t\in\bar I$. Let $\Gamma_t:=\Gamma_{h(t)}$. Let us fix an instant $t\in I$ and consider the space
\begin{equation*}
\begin{aligned}
X_{\Gamma_t}:=\big\{u\in H^1_0(\Omega;\R^3)\,|\,\Div u=0 \text{ in }\Omega\setminus \Gamma_t,\ \DIV u=0 \text{ on } \Gamma_t,\\
\ P_{\Gamma_t}[u]_{\Gamma_t}\in H^1(\Gamma_t;T\Gamma_t)\big\}
\end{aligned}
\end{equation*}
endowed with the canonical scalar product and with the bilinear form
\begin{equation*}
\begin{aligned}
B_{\Gamma_t}(u,\phi)=2\mu_b\int_{\Omega\setminus \Gamma_t} \langle Du,D\phi\rangle_e\,dx+2\mu\int_{\Gamma_t} \langle \D u,\D \phi\rangle_g\,dA.
\end{aligned}
\end{equation*}
We saw in the proof of Theorem 3.11 in \cite{lengeler} that $B_{\Gamma_t}$ defines a uniformly equivalent scalar product on $X_{\Gamma_t}$ for  $h(t)$ being uniformly bounded in $W^{5-4/p}_p(\Gamma)$ with $\|h(t)\|_{L_\infty(\Gamma)}\le\kappa_\Gamma/2$. Also recall from \cite{lengeler} that the weak solution $u\in X_{\Gamma_t}$ of \eqref{eqn:final} is characterized by the equation $B_{\Gamma_t}(u,\phi)=-dF_{\Gamma_t}([\phi]_{\Gamma_t}\cdot\nu_t)$ for all $\phi\in X_{\Gamma_t}$; cf. \eqref{gradient}. Thus, we have
\begin{equation}\label{bgamma}
\begin{aligned}
-\frac{d}{dt}(F(\Gamma_t)-F(\Gamma))&=-dF_{\Gamma_t}([u]_{\Gamma_t}\cdot\nu_t)=B_{\Gamma_t}(u,u)\\&=\|d F_{\Gamma_t}\|_{B'_{\Gamma_t}}^2\ge \frac1c \|dF_{\Gamma_t}\|_{X'_{\Gamma_t}}^2\ge \frac1c \|dF_{\Gamma_t}\|_{Y'_{\Gamma_t}}^2
\end{aligned}
\end{equation}
for some uniform constant $c>0$. Here, the third identity reflects the fact that Riesz' isomorphism is an isometry, for the first estimate we employed the inequality $\|\cdot\|_{B_{\Gamma_t}}\le c \|\cdot\|_{X_{\Gamma_t}}$, and the second estimate follows from the inequality $\|\tilde u\|_{X_{\Gamma_t}}\le c \|\tilde w\|_{Y_{\Gamma_t}}$, where $\tilde u\in X_{\Gamma_t}$ is a suitable extension of $\tilde w\in Y_{\Gamma_t}$, with a uniform constant $c>0$ (this follows from Theorem 3.12 in \cite{lengeler}). Without restriction we can choose $\sigma$ so small that Theorem \ref{thm:ls} can be applied, and thus we infer from \eqref{bgamma} that
\begin{equation*}
\begin{aligned}
\frac{d}{dt}(F(\Gamma_t)-F(\Gamma))\le -\frac1c (F(\Gamma_t)-F(\Gamma))^{2(1-\theta)};
\end{aligned}
\end{equation*}
here, we additionally assumed $\sigma$ to be so small that $F(\Gamma_t)\ge F(\Gamma)$. This shows that for all $t\in I$ we have
\begin{equation}\label{FKonv}
 F(\Gamma_t)-F(\Gamma)\le c\, \left\{\begin{array}{l@{\quad, \quad}l} t^{-1/(1-2\theta)} & \mbox{if }\theta<1/2, \\ e^{-c_0 t} & \mbox{if }\theta=1/2. \end{array}\right.
\end{equation}
Moreover, we can apply Theorem \ref{thm:ls} once more to obtain
\begin{equation}\label{en}
\begin{aligned}
 -\frac{d}{dt}(F(\Gamma_t)-F(\Gamma))^\theta&=\theta (F(\Gamma_t)-F(\Gamma))^{\theta-1}B_{\Gamma_t}(u,u)\\
 &\ge \theta/c\, (F(\Gamma_t)-F(\Gamma))^{\theta-1} \|dF_{\Gamma_t}\|_{Y'_{\Gamma_t}}\|[u]_{\Gamma_t}\cdot\nu_t\|_{Y_{\Gamma_t}}\\
 &\ge \theta/c\, \|[u]_{\Gamma_t}\cdot\nu_t\|_{Y_{\Gamma_t}}.
 \end{aligned}
\end{equation}
Here, for the first estimate we used the identity $B_{\Gamma_t}(u,u)=\|dF_{\Gamma_t}\|_{B'_{\Gamma_t}}\|u\|_{B_{\Gamma_t}}$, the inequalities in \eqref{bgamma}, the uniform equivalence of the scalar products on $X_{\Gamma_t}$, and the uniform continuity of the trace operator $X_{\Gamma_t}\rightarrow Y_{\Gamma_t},\ u\mapsto [u]_{\Gamma_t}\cdot\nu_t$ (see Lemma 3.8 in \cite{lengeler}).
Since $d\Psi_{h(t)}^{-1}([u]_{\Gamma_t}\cdot\nu_t)=\partial_th$, by Remark \ref{gamma} it follows that
\begin{equation*}
\int_0^t\|\partial_th\|_{H^{1/2}(\Gamma)}\,ds\le c \int_0^t\|[u]_{\Gamma_t}\cdot\nu_t\|_{Y_{\Gamma_t}}\,ds\le c\,(F(\Gamma_{h_0})-F(\Gamma))^\theta
\end{equation*}
for all $t\in I$; hence
\begin{equation}\label{abY}
\|h(t)\|_{H^{1/2}(\Gamma)}\le \|h_0\|_{H^{1/2}(\Gamma)}+ c\,(F(\Gamma_{h_0})-F(\Gamma))^\theta.
\end{equation}
Furthermore, for sufficiently small $\sigma$ we deduce from Theorem \ref{thm} that our solution $h$ exists at least as long as $h(t)$ remains in $\bar B_\sigma(0)\subset W^{5-4/p}_p(\Gamma)$, and if the latter is true on a time interval $I'=(0,T')$  then 
\begin{equation}\label{abC}
\|h\|_{C([t',T'], W^{5-4/p+\delta}_p(\Gamma))}\le c
\end{equation}
for arbitrarily small $t'>0$ and some constant $c>0$. Now, let 
\[T^*:=\sup\big\{T'>0\,|\,\|h\|_{C(\bar I',W^{5-4/p}_p(\Gamma))}\le\sigma\big\}.\]
By interpolation and \eqref{abC} there exists a $\theta'\in (0,1)$ such that
\begin{equation*}
\begin{aligned}
\|h(t)\|_{W^{5-4/p}_p(\Gamma)}&\le \|h(t)\|_{W^{5-4/p+\delta}_p(\Gamma)}^{\theta'}\|h(t)\|_{H^{1/2}(\Gamma)}^{1-\theta'}\le c\,\|h(t)\|_{H^{1/2}(\Gamma)}^{1-\theta'}
 \end{aligned}
\end{equation*}
for all $t\in [0,T^*]$. However, from this estimate and \eqref{abY}, for sufficiently small $\epsilon$ we obtain that $\|h(t)\|_{W^{5-4/p}_p(\Gamma)}\le \sigma/2$. Hence, $T^*=\infty$ and the solution exists for all times. Furthermore, from \eqref{abY} we obtain $\partial_t h\in L^1((0,\infty);H^{1/2}(\Gamma))$, and thus $h(t)\rightarrow h_\infty$ in $H^{1/2}(\Gamma))$ for $t\rightarrow\infty$ and some $h_\infty\in H^{1/2}(\Gamma))$. Integrating \eqref{en} over $(t,\infty)$ yields
\begin{equation*}
\begin{aligned}
\|h(t)-h_\infty\|_{ H^{1/2}(\Gamma)}\le c\int_t^\infty\|\partial_sh(s)\|_{ H^{1/2}(\Gamma)}\,ds\le c\,(F(\Gamma_t)-F(\Gamma))^\theta\le c\,t^{-\beta'}
 \end{aligned}
\end{equation*}
for some $\beta'>0$. By weak compactness, we have $h_\infty\in W^{5-4/p+\delta}_p(\Gamma)$ and thus we can use interpolation once more to obtain \eqref{konvergenz}. Finally, from \eqref{konvergenz} and \eqref{FKonv} we conclude that $F(\Gamma_{h_\infty})=F(\Gamma)$.
\qed

\appendix
\section*{Appendix. Second variation of the Canham-Helfrich energy}

We give a brief derivation of the second variation of the Canham-Helfrich energy and of the linearization of its $L_2$-gradient. Of course, this issue has been adressed before in the literature; see, for instance, \cite{guven03,peterson,helfrich2}. Our result agrees with the one in \cite{peterson} (which, however, is not explicitly covariant; see (42), (43)), but it does not agree with the ones in \cite{guven03} (which lacks one of the terms in our expression) and in \cite{helfrich2} (which seems to contain typos; expressions like $Hg^{\alpha\beta}-Kk^{\alpha\beta}$ don't make sense from the point of view of physical dimensions).

Let $\Gamma_t\subset\R^3$ be a closed vesicle moving with velocity $u=w\,\nu_t$. Consider the family of diffeomorphisms $\phi_{t,s}:\Gamma_{t}\rightarrow\Gamma_s$ associated with the vector field $u$, that is, $\phi_{t,t}=\id_{\Gamma_{t}}$ and $\pa_s\phi_{t,s}=u\circ\phi_{t,s}$. We denote by $Df/Dt$ the material derivative of a scalar field $f$ on $\Gamma_t$ with respect to the vector field $u$, that is,
\[\frac{Df}{Dt}\Big|_t:=\pa_s|_{s=t}f\circ\phi_{t,s}.\]
Throughout this appendix we shall work in convected coordinates $(x^\alpha)$ on $\Gamma_t$, that is, $x^\alpha|_t=x^\alpha|_s\circ\phi_{t,s}$. Taking the material derivative of tensor components in such coordinates yields the Lie derivative of the corresponding tensor field which is again a tensor field. In \cite{lengeler} we proved the identities
\begin{equation*}
\begin{aligned}
\frac{D}{Dt}g_{\alpha\beta}=-2 w\, k_{\alpha\beta},&\quad \frac{D}{Dt}g^{\alpha\beta}=-g^{\alpha\mu}g^{\beta\nu}\frac{D}{Dt}g_{\mu\nu}=2 w\, k^{\alpha\beta},\\ \frac{D}{Dt}k_{\alpha\beta}=w_{;\alpha\beta}-w\, k_{\alpha\gamma}k^{\gamma}_\beta,&\quad \frac{D}{Dt}k^{\alpha\beta}=w_;^{\,\alpha\beta}+3w\, k^{\alpha\gamma}k_{\gamma}^\beta,\\
  \frac{D}{Dt}dA=-w\,H\,dA,&\quad \frac{D}{Dt}H=\Delta_g w + w (H^2-2K)
 \end{aligned}
\end{equation*}
The Christoffel symbols are given by
\[\Gamma_{\alpha\beta}^\gamma=\frac12g^{\gamma\delta}(g_{\beta\delta,\alpha}+g_{\alpha\delta,\beta}-g_{\alpha\beta,\delta}).\]
Recall that the difference of two connections is a tensor field and, hence, so is the material derivative of a connection. Thus, we have
\begin{equation*}
\begin{aligned}
\frac{D}{Dt}\Gamma_{\alpha\beta}^\gamma&=-\frac12 \Gamma_{\alpha\beta}^\mu g^{\gamma\nu}\frac{D}{Dt}g_{\mu\nu}+\frac12g^{\gamma\delta}\bigg(\Big(\frac{D}{Dt}g_{\beta\delta}\Big)_{,\alpha}+\Big(\frac{D}{Dt}g_{\alpha\delta}\Big)_{,\beta}-\Big(\frac{D}{Dt}g_{\alpha\beta}\Big)_{,\delta}\bigg)\\
&=\frac12g^{\gamma\delta}\bigg(\Big(\frac{D}{Dt}g_{\beta\delta}\Big)_{;\alpha}+\Big(\frac{D}{Dt}g_{\alpha\delta}\Big)_{;\beta}-\Big(\frac{D}{Dt}g_{\alpha\beta}\Big)_{;\delta}\bigg)\\
&=-k_\alpha^\gamma w_{,\beta}-k_\beta^\gamma w_{,\alpha} + k_{\alpha\beta}w_{,}^{\,\gamma}-w\,k_{\alpha;\beta}^\gamma,
 \end{aligned}
\end{equation*}
where the second identity follows by making use of Riemannian normal coordinates in which the Christoffel symbols vanish (at the center). For a scalar field $f$ on $\Gamma_t$ we have $\Delta_gf=g^{\alpha\beta}(f_{,\alpha\beta} - \Gamma_{\alpha\beta}^\gamma f_{,\gamma})$ and thus
\begin{equation*}
\begin{aligned}
\Big[\frac{D}{Dt},\Delta_g\Big]f&= 2w\,k^{\alpha\beta}f_{;\alpha\beta}-g^{\alpha\beta}f_{,\gamma}\frac{D}{Dt}\Gamma_{\alpha\beta}^\gamma \\
&=2w\,k^{\alpha\beta}f_{;\alpha\beta}+2k^{\alpha\beta}f_{,\alpha}w_{,\beta}-H f_{,\alpha}w_{,}^{\,\alpha}+w f_{,\alpha}H_{,}^{\,\alpha},
\end{aligned}
\end{equation*}
where $[\,\cdot\, ,\cdot\,]$ denotes the commutator. We conclude that
\begin{equation}\label{h1}
\begin{aligned}
\frac{D}{Dt}\Delta_g H&=2w\,k^{\alpha\beta}H_{;\alpha\beta} + 2k^{\alpha\beta}w_{,\alpha}H_{,\beta} - Hw_{,\alpha}H_{,}^{\,\alpha} + w H_{,\alpha}H_{,}^{\,\alpha}\\
&\quad+\Delta_g(\Delta_g w+w(H^2-2K))\\
&=\Delta_g^2w+ \Delta_g w\,(H^2-2K) + w_{,\alpha}\big(2k^{\alpha\beta}H_{,\beta}+3H H_{,}^{\,\alpha}-4K_{,}^{\,\alpha}\big)\\
&\quad + w\big(2k^{\alpha\beta}H_{;\alpha\beta} + H_{,\alpha} H_{,}^{\,\alpha}+\Delta_g(H^2-2K)\big)
\end{aligned}
\end{equation}
Furthermore, we recall that the Riemannian curvature tensor $R$ is given by
\[R_{\alpha\beta\gamma}^{\quad\ \, \delta}=\Gamma_{\gamma\beta,\alpha}^\delta-\Gamma_{\gamma\alpha,\beta}^\delta+\Gamma_{\mu\alpha}^\delta\Gamma_{\gamma\beta}^\mu-\Gamma_{\mu\beta}^\delta\Gamma_{\gamma\alpha}^\mu,\]
while the Ricci tensor $Rc$ satisfies $Rc_{\beta\gamma}=R_{\alpha\beta\gamma}^{\quad\ \, \alpha}=Kg_{\beta\gamma}$. Thus, making again use of Riemannian normal coordinates, we derive
\[\frac{D}{Dt}Rc_{\beta\gamma}=\Big(\frac{D}{Dt}\Gamma_{\gamma\beta}^\delta\Big)_{;\alpha}-\Big(\frac{D}{Dt}\Gamma_{\gamma\alpha}^\delta\Big)_{;\beta},\]
and hence
\begin{equation*}
\begin{aligned}
\frac{D}{Dt}K=\frac{D}{Dt}\frac12\big(Rc_{\beta\gamma}g^{\beta\gamma}\big)&=w\,k^{\beta\gamma}Rc_{\beta\gamma} + \frac12 g^{\beta\gamma}\bigg(\Big(\frac{D}{Dt}\Gamma_{\gamma\beta}^\alpha\Big)_{;\alpha}-\Big(\frac{D}{Dt}\Gamma_{\gamma\alpha}^\alpha\Big)_{;\beta}\bigg)\\
&=w\,K H + H\Delta_g w - k^{\alpha\beta}w_{;\alpha\beta}.
\end{aligned}
\end{equation*}
From this we deduce
\begin{equation}\label{h2}
\begin{aligned}
\frac{D}{Dt}\,H&(H^2/2-2K)\\
&=\big(\Delta_gw+w(H^2-2K)\big)(H^2/2-2K)\\
&\quad+H\big(H(\Delta_gw+w(H^2-2K))-2w\,K H-2H\Delta_g w+ 2k^{\alpha\beta}w_{;\alpha\beta}\big)\\
&=w_{;\alpha\beta}\big((-H^2/2-2K)g^{\alpha\beta}+2Hk^{\alpha\beta}\big)+w(3H^4/2-7KH^2+4K^2)
\end{aligned}
\end{equation}
and 
\begin{equation}\label{h3}
\begin{aligned}
\frac{D}{Dt}\,(2K-HC_0/2)&=w_{;\alpha\beta}\big((2H-C_0/2)g^{\alpha\beta}-2k^{\alpha\beta}\big)\\
&\quad + w(2KH-C_0/2\,H^2+C_0K).
\end{aligned}
\end{equation}
Collecting \eqref{h1}, \eqref{h2}, and \eqref{h3}, we finally compute the linearization of $\grad_{L_2}F_{\Gamma_t}$ by taking its material derivative, that is, 
\begin{equation*}
\begin{aligned}
\frac1\kappa\frac{D}{Dt}&\grad_{L_2}F_{\Gamma_t}\\
&=\Delta_g^2w+w_{;\alpha\beta}\big((H^2/2-4K+2HC_0-C_0^2/2)g^{\alpha\beta}+2(H-C_0)k^{\alpha\beta}\big)\\
&\quad + w_{,\alpha}\big(2k^{\alpha\beta}H_{,\beta} + 3HH_{,}^{\,\alpha}-4K_{,}^{\,\alpha}\big)+ w\big(2k^{\alpha\beta}H_{;\alpha\beta} + \Delta_g(H^2-2K)\\
&\quad+ H_{,\alpha}H_,^{\,\alpha}+3H^4/2 - 7KH^2+4K^2 +2C_0KH-C_0^2/2\,H^2+C_0^2K\big)\\
&=:\Delta_g^2w+(a^{\alpha\beta}w_{,\alpha})_{;\beta} + \tilde b\,w.
\end{aligned}
\end{equation*}
Since
\[\frac{d^2}{dt^2}F(\Gamma_t)=d^2F_{\Gamma_t}(w,w)+dF_{\Gamma_t}\Big(\frac{D}{Dt}w\Big),\]
the second variation of the Canham-Helfrich energy is given by
\begin{equation*}
\begin{aligned}
d^2F_{\Gamma_t}(w,w)&=\int_{\Gamma_t}\frac{D}{Dt}\grad_{L_2}F_{\Gamma_t}\ w\,dA+\int_{\Gamma_t}\grad_{L_2}F_{\Gamma_t}\,w\,\frac{D}{Dt}dA\\
&=\kappa\int_{\Gamma_t}\big((\Delta_gw)^2+a^{\alpha\beta}w_{,\alpha}w_{,\beta} + b\,w^2\big)\,dA,
\end{aligned}
\end{equation*}
where \[b=((2k^{\alpha\beta}-Hg^{\alpha\beta})H_{,\alpha})_{;\beta} +\Delta_g(H^2-2K)+H^4-5KH^2+4K^2 +C_0^2K;\]
note that the first term on the right hand side is not contained in the analogous expression (84) in \cite{guven03}.

\bibliographystyle{plain}
\bibliography{references}
\end{document}